\documentclass{amsart}
\usepackage{amsmath,amscd,amssymb,amsthm,epsfig,psfrag,latexsym,subfigure}

\newtheorem{theorem}{Theorem}[section]
\newtheorem{lemma}[theorem]{Lemma}

\newtheorem{definition}[theorem]{Definition}

\newtheorem*{acknowledgments}{Acknowledgments}

\newcommand{\NN}{\mathbb{N}}

\newcommand{\CC}{\mathbb{C}}
\newcommand{\RR}{\mathbb{R}}

\def\vol{\mbox{\rm{Vol}}}
\def\length{\mbox{\rm{length}}}
\def\area{\mbox{\rm{Area}}}

\def\a{\alpha}
\def\b{\beta}
\def\d{\partial}
\def\D{\Delta}
\def\e{\epsilon}
\def\g{\gamma}
\def\i{\iota}

\def\l{\lambda}
\def\k{\kappa}

\def\i{\iota}

\def\S{\Sigma}

\def\tr{\mbox{\rm{tr}}}

\def\Int{\mbox{int}\ }
\def\exp{\mbox{\rm{exp}}\ }

\def\HH{\mathbb{H}}

\def\MN{\mathcal{N}}
\def\MO{\mathcal{O}}

\def\CC{\mathbb{C}}
\def\CH{\mathcal{CH}}

\edef\t@mp{\catcode`\noexpand\#=\the\catcode`\#}%
    \catcode`\#=12
    \def\h@sh{#}%
    \t@mp

\edef\t@mp{\catcode`\noexpand\~=\the\catcode`\~}%
    \catcode`\~=12
    \def\tild@{~}%
    \t@mp

\begin{document}
\Large
\title{Tameness of hyperbolic 3-manifolds}
\author{Ian Agol}
\address{Department of Mathematics, University of Illinois at
Chicago, 322 SEO m/c 249, 851 S. Morgan Street, Chicago, IL
60607-7045}

\email{agol@math.uic.edu}
\thanks{Partially supported by ARC grant 420998, NSF grant DMS 0204142, and the Sloan
Foundation}
\maketitle
\section{Introduction}
Marden conjectured that a hyperbolic 3-manifold $M$ with finitely
generated fundamental group is tame, {\it i.e.} it is homeomorphic
to the interior of a compact manifold with boundary \cite{Marden}.
Since then, many consequences of this conjecture have been
developed by Kleinian group theorists and 3-manifold topologists.
We prove this conjecture in theorem \ref{tame}, actually in
slightly more generality for PNC manifolds with hyperbolic cusps
(the cusped case is reduced to the non-cusped case in section
\ref{cusps}, see the next section for definitions of this
terminology).

Many special cases of Marden's conjecture have been resolved and
various criteria for tameness have been developed, see
\cite{Marden, Th, Bon, Ca94, CM96, FM98,  Souto02, BBES, BS04,
KS03}. Bonahon resolved the case where $\pi_1(M)$ is
indecomposable, that is $\pi_1(M)\neq A*B$. This was generalized
by Canary \cite{Ca89, Ca93} to the case of PNC manifolds with
indecomposable fundamental group. Various other cases of limits of
tame manifolds being tame have been resolved, culminating in the
proof by Brock and Souto that algebraic limits of tame manifolds
are tame \cite{ Th, CM96, BBES, BS04}. In a certain sense,
Canary's covering theorem provides other examples of tame
manifolds \cite{Ca}, whereby tameness of a cover implies tameness
of a quotient under certain circumstances. Conversely, Canary and
Thurston showed that covers (with finitely generated fundamental
group) of tame hyperbolic manifolds of infinite volume are tame
(this is a purely topological result) \cite{Ca94}.

Bonahon in fact proved that hyperbolic 3-manifolds with
indecomposable fundamental group are geometrically tame
\cite{Bon}. Geometric tameness was a condition formulated by
Thurston (and proved in some special cases by him) which implies
that the ends are either geometrically finite or simply degenerate
\cite{Th}. Canary showed that tame manifolds are geometrically
tame \cite{Ca93}. He also showed (generalizing an argument of
Thurston \cite{Th}) that geometric tameness implies the Ahlfors
measure conjecture \cite{Ahlfors64}. Canary's arguments were
generalized for PNC manifolds by Yong Hou \cite{Hou03}. In fact,
these arguments give a geometric proof of the Ahlfors finiteness
theorem \cite{Ahlfors64}. Other corollaries for Kleinian groups
are noted in \cite{Ca94}.

Assuming the geometric tameness conjecture, Thurston conjectured
that a hyperbolic 3-manifold $N$ is determined by its end
invariants: the topological type of $N$, the conformal structure
of the domain of discontinuity of  each geometrically finite end,
and the ending lamination of each simply degenerate end. This was
resolved last year for tame manifolds by Brock, Canary and Minsky
(with major contributions from Masur-Minsky) \cite{Minsky01,
Minsky03, Minsky, MaMi99, MaMi00, Minsky94, BCM}. This then gives
a complete parameterization of isometry types of Kleinian groups.
A corollary of the ending lamination conjecture is the density
conjecture: any Kleinian group is an algebraic limit of
geometrically finite Kleinian groups (shown to be a corollary of
the ending lamination conjecture by Kleineidam-Souto \cite{KS02}).

An application of the covering theorem and geometric tameness
(which, as noted, follows from tameness) implies that a cover with
finitely generated fundamental group of a finite volume hyperbolic
3-manifold is either geometrically finite or is associated to the
fiber of a fibration over $S^1$ (Conjecture C \cite{Ca94}). This
has some corollaries. The Simon conjecture \cite{Simon76} says
that irreducible (finitely generated) covering spaces of compact
3-manifolds are tame. It follows from Simon's work and tameness of
hyperbolic 3-manifolds that Simon's conjecture holds for manifolds
satisfying the geometrization conjecture.

 Another application is that the
fundamental group of the figure 8 knot complement is LERF (as well
as for the fundamental groups of several other compact 3-manifold
groups) \cite{ALR, Gitik03, Wise98, LR01, LongReid}.

 It is shown
that the minimal volume orientable hyperbolic 3-manifold has a
finite index subgroup which is generated by two elements in
\cite{CHS}. It is conjectured that the Weeks manifold (which has
volume $=.9427...$ ) is the minimal volume orientable hyperbolic
3-manifold. It is known to be minimal volume among arithmetic
manifolds \cite{CFJR}. Thus, the commensurator of the fundamental
group of the minimal volume orientable hyperbolic 3-manifold must
have torsion, so it must cover a non-trivial orbifold with torsion
(since a 2-generator group has a hyperelliptic involution which
must lie in the commensurator).

Now we outline the argument. The cusped case is reduced to the
cusp-free case in section \ref{cusps}. In section \ref{geometric},
it is shown that a PNC manifold with finitely generated
fundamental group, satisfying an additional technical condition,
has a sequence of homotopy equivalent tame PNC branched covers
which limit to it geometrically. In section \ref{limits}, it is
shown how to conclude in this case that the manifold is tame by
generalizing an argument of Canary-Minsky \cite{CM96}. The general
case is reduced to this case via section \ref{sequence} and the
PNC covering theorem \ref{covering}. We add an appendix which
gives details of an argument of Kleiner \cite{Kleiner92} which we
adapt to show that the boundary of the convex hull of a PNC
manifold has area bounded solely in terms of its Euler
characteristic and the pinching constants (this is used in the
limit argument in section \ref{limits}).

In several theorems, questions about hyperbolic manifolds have
been reduced to questions about Haken manifolds, by ``drilling''
(or branching) along an appropriately chosen geodesic (see
\cite{Ag0},\cite{Ca93},\cite{Gabai94}). A reminiscent procedure of
grafting allows one to construct geometrically finite manifolds
from geometrically infinite ones \cite{BrockBromberg02}. The
arguments in this paper are inspired by these constructions, and
are direct  generalizations of Canary's trick \cite{Ca93} (see
also \cite{Souto02}).

\begin{acknowledgments}
We thank Ken Bromberg, Dick Canary, Bruce Kleiner, Yair Minsky,
Bobby Myers, Manuel Ritor\'e, and Juan Souto for helpful
conversations. Part of the line of argument in section
\ref{limits} was suggested by Minsky. We thank Bobby Myers for
writing up a slight generalization of his result which was
important for our argument. We also thank Mike Freedman, who had
much influence on the arguments in this paper.
\end{acknowledgments}

\section{Definitions}

Let $\HH^3_b$ be hyperbolic space endowed with a Riemannian metric
of constant sectional curvature $b<0$.

If $C\subset M$ is a subset of a metric space, then we'll use
$\MN_R(C)$ to denote the set of points of distance $< R$ to $C$.

A {\it PNC (pinched negatively curved) } manifold is a Riemannian
manifold $(M,g)$, where $g$ is a complete Riemannian metric, and
there exist constants $b<a<0$ such that $b<K(P)<a$, where $P$ is
any 2-plane in $G_2(T M)$ and $K(P)$ is the sectional curvature of
the metric along $P$ (see \cite{BGS85} for more about PNC
manifolds).

Let $\e$ be a Margulis constant for $(M,g)$, $M_{thin(\e)}$ the
set of points at which the injectivity radius is $< \e$ (see 8.3,
10.3 \cite{BGS85} ). A {\it cusp} $Q$ of $M$ is a component of
$M_{thin(\e)}$ such that elements of $\pi_1(Q)$ have arbitrarily
short loop representatives (equivalently, a cusp is a non-compact
component of $M_{thin(\e)}$). We will denote
$M_{thick(\e)}=M_{\e}=M-\Int M_{thin(\e)}$.

A PNC manifold $(M^3,g)$ has {\it hyperbolic cusps} if for $\e$
small enough, every cusp $Q$ of $M$ is locally isometric to
$\HH^3_b$, for some $b<0$. Equivalently $Q$ is isometric to a
quotient of a horoball in $\HH^3_b$ by a discrete group of
parabolic isometries. We will also consider PNC orbifolds, which
have charts locally modelled on a ball in a PNC manifold quotient
a finite group of isometries.

A subset $U\subset M$ is convex if for every geodesic $p$ with
endpoints $x,y\in U$, $p\subset U$. This is not the usual
definition, since we don't require that $p$ be length-minimizing.
Clearly, intersections of convex subsets are convex, so there
exists a smallest convex subset. If $M$ is PNC, then let $\CH(M)$
denote the convex core, that is the smallest closed convex subset
of $M$. $\CH(M)$ is homotopy equivalent to $M$, and is generally a
codimension zero submanifold, unless $M$ has a totally geodesic
invariant subspace.

Let $S$ be a closed surface, with triangulation $\tau$. Let $N$ be
a PNC manifold, and $f:S\to N$ a map. Then $f$ is a {\it
simplicial ruled surface} if $f$ takes each edge of $\tau$ to a
geodesic segment, each face of $\tau$ has a foliation by geodesic
segments (in fact, this foliation is irrelevant for the arguments,
since any such triangular disk will lie in a $\delta$ neighborhood
of the edges of the triangle), and the cone angle around every
vertex of $\tau$ is $\geq 2\pi$ (see \cite{Ca93, Souto02}).

A collection of pairwise non-conjugate primitive elements
$\a\subset \pi_1(M)-\{1\}$ is {\it algebraically diskbusting} if
for any free product representation $\pi_1(M)=A*B$, $A\neq 1 \neq
B$, $\a \nsubseteq A \cup B$. This terminology is used, since if
$C\subset M$ is a compact core, and $\a^*\subset C$ is a link
whose components represent the conjugacy classes of $\a$, then $\d
C$ is incompressible in $M-\a^*$, that is $\a^*$ is diskbusting in
$C$.

Let $M$ be an irreducible, connected, orientable 3-manifold. A
3-manifold with boundary $J\subset M$ is {\it regular} if $J$ is
compact and connected, and $M-\Int J$ is irreducible and has no
compact component.

A {\it regular exhaustion} $\{C_n\}$ of $M$ is a sequence of
regular 3-manifolds such that $C_n\subset \Int C_{n+1}$ and
$M=\cup_n C_n$.

An open submanifold (but not necessarily properly embedded)
$V\supset J$ is an {\it end reduction of M at J} if there exists a
regular exhaustion $\{M_n\}$ with $V=\cup_n M_n$, such that $\d
M_n$ is incompressible in $M-J$, $V-J$ is irreducible, $M-V$ has
no compact components, and any submanifold $N\subset M$ with
$J\subset \Int N$, $\d N$ incompressible in $M-J$ has the property
that $N$ may be isotoped fixing $J$ so that $N\subset V$ (this is
called the {\it engulfing property}).

We need the following generalization of strong convergence of
Kleinian groups, which is a type of algebraic and geometric limit
(the notion of geometric convergence which we use is much stronger
than usual, but suffices for our purposes).
\begin{definition} \label{strong limit}
Let $N_i$ be a sequence of uniformly PNC manifolds (with finitely
generated fundamental group). $N_i$ converges {\it super strongly}
to $M$ if there is a compact core $C\subset M$, and for all
compact connected submanifolds $K$ with $C\subset K \subset M$,
there is an $I$ such that for $i>I$, there is a Riemannian
isometry $\eta_{i,K} :K \to N_i$, with $\eta_{i,K|C}: C\to N_i$ is
a homotopy equivalence.
\end{definition}

\begin{definition} \label{drilling} Let $L\subset M$ be a link in $M$. Then $M_L(2,0)$ is the
meridianal $\pi$-orbifold drilling on $L$, which may be defined in
the following way. Take a regular neighborhood $\MN(L)$, then each
component is a solid torus (assuming $M$ is orientable). If one
has a solid torus $\D=D^2\times S^1$, let $K={0}\times S^1$, the
core of $\D$. There is an involution $\i:D^2\times S^1\to
D^2\times S^1$ given by $\i(x,y)=(-x,y)$, which is $\pi$ rotation
of $D^2 \subset \CC$ and fixes $K$.  Then $\D/\i$ has a natural
orbifold structure , called $\D_K(2,0)$. The underlying space is
still a solid torus, and we may remove each component of $\MN(L)$,
and glue in a copy of $\D_K(2,0)$ so that the meridian orbidisk of
$\D_K(2,0)$ is glued to the boundary of the meridian disk of
$\MN(L)$. This corresponds to gluing each copy of  $\d \D_K(2,0)$
to $\d M-\Int \MN(L) $ by $(x,y)\equiv (x^2,y)$.
\end{definition}

\section{PNC orbifolds}

We record the following lemma, whose proof is exactly like that of
the manifold case, since the exponential map is an orbifold
covering map.
\begin{lemma} (Ghys-Haefliger \cite{GH90})
If $\MO^n$ is a PNC orbifold, then $\widetilde{\MO}$ exists and is
a homeomorphic to $\RR^n$ with an induced PNC metric.
\end{lemma}

We need to use the following branched cover trick, which will
occur several times in the argument.

\begin{lemma} \label{branched lifting}
Let $M^3$ be a PNC manifold, with geodesic link $\a^*$, which has
a neighborhood $\MN_R(\a^*)$ locally isometric to $\HH^3_{b}$,
$b<0$. Let $M_{\a^*}(2,0)$ be the PNC $\pi$-orbifold obtained by
modifying the metric in $\MN_R(\a^*)$ to get a $\pi$-orbifold
along the link $\a^*$ \cite{GrTh87}. Let $\MN_R(\a^*)\subset
U\subset M$ be a submanifold, and let $V\to M_{\a^*}(2,0)$ be an
orbifold cover of $M_{\a^*}(2,0)$. Suppose that $U_{\a^*}(2,0)$
lifts to $V$. Let $V'=V-\MN_R(\a^*)(2,0)\cup \MN_R(\a^*)$ be the
orbifold obtained by removing the orbifold locus and replacing the
original metric on the lift of $\a^*$ to $V$. Then $U$ lifts to
$V'$.
\end{lemma}
\begin{proof}
It's clear that when we remove $\MN_R(\a^*)(2,0)$ from $V$, and
glue back in $\MN_R(\a^*)$ to get $V'$, the same operation occurs
on the suborbifold $U_{\a^*}(2,0)$ to obtain an embedding
$U\subset V'$.
\end{proof}

\section{End reductions}

\begin{theorem} (Brin-Thickstun, 2.1-2.3 \cite{BT87}) \label{end reduction}
 Given $J$ a regular submanifold of $M$, an
end reduction $V$ at $J$ exists and is unique up to non-ambient
isotopy fixing $J$.
\end{theorem}

Let $\a^*\subset M$ be an algebraically diskbusting link, and
$V\subset M$ an end reduction of $M$ at $\MN(\a^*)$.

\begin{theorem} (Myers, 9.2 \cite{Myers04}) \label{Myers}
If $\a^*\subset M$ is an algebraically diskbusting link, and $V$
is an end reduction at $\MN(\a^*)$ with $\i:V\to M$ the inclusion
map, then $\iota_*:\pi_1(V)\to \pi_1(M)$ is an isomorphism. It is
also clear that $\pi_1(V-\a^*)$ injects into $\pi_1(M-\a^*)$.
\end{theorem}

\section{Pared manifold compression bodies}
In \cite{Bonahon83}, Bonahon shows that a compact, irreducible
3-manifold can be decomposed along incompressible surfaces into
finitely many compression bodies and manifolds with incompressible
boundary, which are unique up to isotopy. Thus, if one has a PNC
3-manifold $M$ with finitely generated fundamental group and core
$C$, the ends of $C$ will correspond to boundary components of
$M$, which are either incompressible or compressible. The ends
corresponding to incompressible boundary will be tame, by Canary
\cite{Ca89}. Each end corresponding to  a compressible boundary
component will be associated to a compression body $B$ in
Bonahon's decomposition. Passing to the cover $N$ of $M$
corresponding to the fundamental group of this compression body
$B$, $B$ lifts to a core of $N$, and thus the end will be
associated to the unique compressible boundary component of $B$.
Again, by Bonahon's result \cite{Bon}, the other ends of $N$
corresponding to incompressible boundary components of $B$ will be
tame. Thus, to prove that $M$ is tame, it suffices to prove that
the end of $N$ corresponding to the compressible boundary
component of $B$ is tame.

We assume that $M$ has a complete PNC Riemannian metric with
hyperbolic cusps, whose union we will call $H$. If we choose a
small enough Margulis constant $\e$, then we may assume that $H$
is the union of non-compact components of $M_{thin(\e)}$. Then we
will denote $M-\Int H=M_H$, the neutered manifold. $\d M_H$
consists of flat open annuli or tori locally isometric to a
horosphere in $\HH^3_{-b^2}$. By the relative core theorem
\cite{McC86}, there is a core $(C_H,P)\subset (M_H,\d M_H)$ which
is a homotopy equivalence of pairs, such that $C_H\cap \d M_H=P$.
The manifold $(C_H,P)$ is a pared manifold: $P$ consists of
incompressible annuli and tori in $\d C_H$ such that no two annuli
are parallel in $C_H$ and every torus component of $\d C_H$ is
$\subset P$ (see Def. 4.8, \cite{Mo}).

\section{Reduction of the cusped case to the case with no cusps}
\label{cusps}
 As in the previous section, we assume that the
compact core is a compression body. We'll also assume that the
manifold is orientable, by passing to a 2-fold orientable cover.
For each rank 2 cusp, we perform a high Dehn filling to obtain a
PNC manifold with the same (non-cusp) ends, using the method of
the Gromov-Thurston $2\pi$ theorem (\cite{GrTh87, BH}). The Scott
core will remain a compression body.

So we may assume that all of the cusps are rank 1. For each rank
one cusp $Q$, $\d Q$ is an infinite annulus. Fill along the core
curve of this annulus by a $2\pi/n$ orbifold, where $n$ is large
enough to obtain a PNC orbifold $M'$, again using the method of
the $2\pi$ theorem (\cite{GrTh87, BH}). The Scott core $C'$ of
$M'$ is obtained by adding orbifold 2-handles to the components of
the pared locus $P\subset C_H$. Since $M'$ is PNC, $C'$ will be an
atoroidal, irreducible orbifold. By the orbifold theorem
\cite{BLP, CHK}, $\Int C'$ has a hyperbolic structure. Thus, we
may pass to a finite sheeted manifold cover of $M'$ by Selberg's
lemma. If we can show that this cover is tame, then $M'$ will be
tame by the finite covering theorem \ref{finitecovering}, and thus
our original manifold $M$ is tame .

\section{Metric modulo $\e$-thin part}

Assume that $M^3$ is a PNC manifold with no cusps. We may choose
$\e$ small enough such that each pair of components of
$M_{thin(\e)}$ has distance $> \e$. Define a distance function
$d_{\e}$ on $M_{\e}=M_{thick(\e)}$ by $d_{\e}(x,y)
=\underset{p}{\inf}\ \length(p\cap M_{\e})$, where $p$ is a path
connecting $x$ and $y$. $d_{\e}$ is a metric on $\Int M_{\e}$,
whose completion is obtained from $M$ by crushing components of
$M_{thin(\e)}$ to points. Let $\MN_R^{\e}(C)$ denote the $R$
neighborhood of $C$ in the metric $d_{\e}$.

\begin{lemma} \label{diameter}
If $S\subset M$, and $diam(S)=\infty$, then $diam_{d_{\e}}(S\cap M_{\e})=\infty$.
\end{lemma}
\begin{proof}
Fix $x\in S \cap M_{\e}$. Then if $d_{\e}(x,y)\leq D$, and $p$ is
a path in $M$ joining $x$ and $y$ such that $length(p\cap
M_{\e})=d_{\e}(x,y)$, then $p$ may intersect at most $D/\e$
components of $M_{thin(\e)}$, since the distance between
components is at least $\e$. Also, it is clear that each segment
of $p\cap M_{\e}$ has $\length\leq D$. Let $N=\lfloor D/\e
\rfloor$, and inductively define for $0\leq k\leq N$, $B_0=\{x\}$,
and $B_k=\MN_D(B_{k-1} \cup U_k)$, where $U_k$ is the collection
of components of $M_{thin(\e)}$ meeting $B_{k-1}$. If $B_{k-1}$ is
compact, then $U_k$ has finitely many components, so $B_k$ is also
compact. Thus, $B_N$ is compact, therefore has bounded diameter,
and clearly $\MN^{\e}_{d_{\e}(x,D)}(x)\subset B_N$, so we see that
$\MN^{\e}_{d_{\e}(x,D)}(x)$  has bounded diameter. Thus, $S \cap
M_{\e} \nsubseteq \MN^{\e}_{d_{\e}(x,D)}(x)$ for any $D$, so
$diam_{d_{\e}}(S \cap M_{\e})=\infty$.
\end{proof}

\section{Bounded diameter lemma}
The following lemma in various forms is due to Thurston (8.8.5
\cite{Th}), Bonahon (1.10 \cite{Bon}), Canary (3.2.6 \cite{Ca89}),
and Canary-Minsky (5.1 \cite{CM96}). We state it for PNC orbifolds
$M$ with no cusps. The proof is essentially the same as the
previous proofs.

\begin{lemma}
Given $A\in \NN$, and a Margulis constant $\e$, there are
constants $D$ and $\mu$ satisfying the following: if $f:\S\to M$
is a  simplicial ruled surface with $-\chi(\S) \leq A$, such that
every compressible or accidental elliptic curve has length $>\mu$,
then $diam_{d_{\e}}(f(\S)\cap M_{\e}) \leq D$.
\end{lemma}

\section{Drilling and geometric limits} \label{geometric}
 Let $(M,\nu)$ be a connected orientable 3-manifold with
a PNC metric $\nu$ with no cusps such that $\pi_1(M)$ is finitely
generated. We may assume that $M$ has a compact core $C$ which is
a compression body. Let $\a\subset \pi_1(M)$ be a collection of
algebraically diskbusting elements. Let $\a^*$ be a union of
geodesic representatives of (the conjugacy classes of) $\a$ in $M$
(we know that each curve will have distinct representatives, since
the elements of $\a$ are primitive and non-conjugate). Perturb the
metric on $M$ near $\a^*$ to a PNC metric $\nu'$ such that $\a^*$
is homotopic to an embedded geodesic $\a'$, by the method of Lemma
5.5, \cite{Ca93}. Then we may deform the metric near $\a'$ to a
metric $\nu''$ so that the metric near $\a''$ (the geodesic
representative of $\a$ in $\nu''$) is locally isometric to
$\HH^3_b$, $b<0$, by Lemma 3.2 and 4.2, \cite{Hou03}. To simplify
notation, we will denote this new metric by $\nu$, and we will
assume that the geodesic representative of $\a$ is a geodesic link
$\a^*$ with neighborhood locally isometric to $\HH^3_b$.

Let $V$ be an end reduction of $M$ at $\MN(\a^*)$, let $\{M_i\}$
be a regular exhaustion of $V$, $\a^*\subset M_0$, so that $\d
M_i$ is incompressible in $M-\a^*$. We may choose a compact core
$C\subset V$, and assume that $C\subset M_0$.

\begin{definition} Let $V\subset M$ be an open irreducible submanifold such that
$V=\cup_i M_i$ is an exhaustion by regular submanifolds containing
$\MN(\a^*)$ and so that inclusion $V\hookrightarrow M$ induces an
isomorphism $\pi_1 V\to \pi_1 M$, and $\d M_i$ is incompressible
in $\pi_1(M-\a^*)$. We call such a submanifold an {\it almost end
reduction}.
\end{definition}
These conditions follow if $V$ is an
end reduction of $M$ at $\a^*$ by \ref{Myers}, since $\a^*$ is
algebraically diskbusting, thus an end reduction is an almost end
reduction in this case.

\begin{theorem} \label{case1}
Assume that $\pi_1(V_{\a^*}(2,0))\to\pi_1(M_{\a^*}(2,0))$ is an
isomorphism, where $V$ is an almost end-reduction (recall def.
\ref{drilling} for the notation).
 Under the above hypotheses, $M$ is a super strong limit of PNC tame
 $N_i$.
\end{theorem}
\begin{proof}
For simplicity, let's first consider the case that $M$ can be
exhausted by compact cores $M=\cup_{i} C_i$ (in which case, we may
assume that $V=M$). This case was shown by Souto to be tame
\cite{Souto02}, but we will outline how our approach works in this
case, since several aspects of the argument are simplified (this
argument is actually closely related to Souto's).

The orbifold $M_{\a^*}(2,0)$ is exhausted by suborbifolds
$C_{i,\a^*}(2,0)$ (for $i\gg 0$) such that
$\pi_1(C_{i,\a^*}(2,0))\hookrightarrow \pi_1(M_{\a^*}(2,0))$.
Since $\pi_1(C_i)\simeq \pi_1(M)$, $\partial C_i$ is
incompressible in $M-\Int(C_i)$. If $\pi_1(\partial C_i)$ does not
inject into $C_{i,\a^*}(2,0)$, then by the equivariant Dehn's
lemma, there is an orbifold disk $(D,\d D)\subset
(C_{i,\a^*}(2,0),\d C_i)$. $D$ cannot be disjoint from $\a^*$,
since $\a$ is algebraically diskbusting. So $D$ must be an
orbifold disk, which meets $\a^*$ exactly once. But this means
that $\a^*$ represents a generator of $\pi_1(C_i)$, which means
that it is not algebraically diskbusting, a contradiction.

Now, let $C'_i \to M_{\a^*}(2,0)$ be the orbifold cover such that
$\pi_1(C'_i)=\pi_1(C_{i,\a^*}(2,0))$. Then $\pi_1(C'_i)$ is
indecomposable, so by Canary's theorem \cite{Ca93}, $C'_i$ is
tame. We have a lift $C_{i,\a^*}(2,0)\hookrightarrow C'_i$, which
is a homotopy equivalence of orbifolds. We may create a branched
cover of $M$ by taking $N_i=C'_i-C_{i,\a^*}(2,0) \cup C_i$, where
we glue along $\d C_i$ (this corresponds to ``erasing'' the
orbifold locus of $C'_i$, as in lemma \ref{branched lifting}).
$N_i$ has a PNC metric, since we have really only inserted the
original metric in $\MN(\a^*)$ back into $C'_i$.

The claim is that the manifolds $N_i$ limit super strongly to $M$
(recall def. \ref{strong limit}). To see this, let $K\subset M$ be
a compact subset. Then we may choose $i$ large enough that
$K\subset C_i$. Then $C_i$ lifts isometrically up to $N_i$ by
lemma \ref{branched lifting}, and therefore so does $K$. If
$C\subset K$ is  a compact core, then the lift will induce a
homotopy equivalence $\pi_1(C)\simeq \pi_1(C_i)\simeq \pi_1(N_i)$.

In general $M$ will not be exhausted by compact cores. To try to
mimic the above proof, we use  the notion of the almost end
reduction.

\begin{lemma}
 $\pi_1(M_{i,\a^*}(2,0))\to
\pi_1(M_{\a^*}(2,0))$ is an injection for large enough $i$.
\end{lemma}
\begin{proof}
By the definition of almost end reduction, $\d M_i$ is
incompressible in $M-\Int(M_i)$ and in $M-\a^*$. If $\pi_1(\d
M_i)$ does not inject into $\pi_1(M_{i,\a^*}(2,0))$, then by the
equivariant Dehn's lemma \cite{MeeksYau81}, there exists an
orbifold disk $D_i$ which meets $\a^*$ exactly once (since $M_i$
is a compression body, we may take the 2-fold branched cover over
$\a^*$. Since $\a^*$ does not meet $\d M_i$, the only possibility
is for the 2-fold action to fix a disk by rotation about a point.
One may also prove this using the orbifold theorem \cite{BLP,
CHK}). In $M_i'=M_i-\Int \MN(\a^*)$, we get an annulus
$A_i=D_i\cap M_i'$ such that one boundary component lies on $\d
M_i$ and one lies on a meridian of $\d \MN(\a^*)$. Then for $j>i$,
$A_j$ may be isotoped so that $A_j\cap M_i'$ is a collection of
essential annuli, and has exactly one annulus with a boundary
component a meridian of $\d \MN(\a^*)$. Let $H_i$ be the
characteristic product region of $M_i'$ \cite{JS79, Johannson79},
then we may assume that $A_i\subset H_i$. There are finitely many
components of $H_i$ of the form $S^1\times [0,1]^2$ which meet $\d
\MN(\a^*)$ (exercise), and any annulus with one boundary on a
meridian of $\d \MN(\a^*)$ and the other on $\d M_i$ must
intersect the core of one of these annuli. We may take a subset
$W\subset \NN$ such that for $i\in W$, for all $j\in W, j>i$,
$A_j\cap M_i'$ has the same annulus meeting $\d \MN(\a^*)$ up to
isotopy (which we relabel to be $A_i$, and relabel $W$ to be
$\NN$), and such that $(|A_j\cap M_{j-1}'|,...,|A_j\cap
M_1'|,|A_j\cap M_0'|)$ is minimal with respect to lexicographic
order.

By Haken finiteness, there are only finitely many disjoint
non-parallel annuli in $A_j\cap M_i'$. We may assume that no two
are parallel, otherwise we could cut and paste to get an annulus
$A_j'$ in $M_j'$ meeting $M_i'$ in fewer components, which
decreases $(|A_j\cap M_{j-1}'|,...,|A_j\cap M_1'|,|A_j\cap M_0'|)$
with respect to lexicographic order. Since for $k>j>i$, $A_k\cap
M_j' \supset A_j$, we see that $A_j\cap M_i'$ is isotopic to a
subset of $A_k\cap M_i'$. Thus, $A_j\cap M_i'$ stabilizes. Then we
may form the union $\cup_i A_i$ to get a properly embedded
punctured disk in $V-\Int \MN(\a^*)$, which capping off with a
meridian disk of $\MN(\a^*)$, gives a properly embedded plane
$P\subset V$ meeting $\a^*$ exactly once. But this implies that
the component $\a'$ of $\a^*$ meeting $P$ is a free generator of
$\pi_1(V)$, and therefore $\a$ is not a collection of
algebraically disk-busting elements, since $\pi_1(V)=\langle \a'
\rangle * \pi_1(V-P)$, where $\a-\a'\subset \pi_1(V-P)$.

\end{proof}

Given this lemma, we may assume that $C\subset M_0$ and
$\pi_1(M_{i,\a^*}(2,0)) \hookrightarrow \pi_1(M_{\a^*}(2,0))$. Let
$M_i'\to M_{\a^*}(2,0)$ be the orbifold covering such that
$\pi_1(M_i')= \pi_1(M_{i,\a^*}(2,0))$, with corresponding lift
$M_{i,\a^*}(2,0) \hookrightarrow M_i'$. Let $M_i''$ be the
manifold which is obtained from $M_i'$ by deleting the orbifold
locus. That is, we replace $\MN(\a^*)(2,0)\subset M_i'$ by
$\MN(\a^*)$ with its original metric, to get the PNC manifold
$M_i''$. Then $M_i''$ is a branched cover of $M$, branched over
$\a^*$. We have a lift $M_i\hookrightarrow M_i''$ obtained by
extending $M_i-\MN(\a^*) \hookrightarrow M_i'$ over $\MN(\a^*)$,
by lemma \ref{branched lifting}. Since $C\subset M_0\subset M_i$,
we have a lift $C\hookrightarrow M_i''$, so let $N_i$ be the cover
of $M_i''$ corresponding to $im\{\pi_1(C)\to \pi_1(M_i'')\}$. Note
that this construction will be trivial if $V$ is itself tame, and
corresponds to the same construction as described above when $M$
is exhausted by compact cores.

Let $K\subset M$ be a compact connected submanifold, and assume
$\MN(\a^*) \subset \Int K$. Since $C$ is a compact core for $M$,
we may homotope $K$ into $C$, such that the homotopy lies in a
connected compact submanifold $K'\subset M$ (where we assume
$C\subset K'$). Thus, $im\{\pi_1(K)\to \pi_1(K')\} \subset
im\{\pi_1(C)\to \pi_1(K')\}$. Since we are assuming that
$\pi_1(M_{\a^*}(2,0))=\pi_1(V_{\a^*}(2,0))=\cup_i
\pi_1(M_{i,\a^*}(2,0))$, there exists $i$ such that
$\pi_1(K'_{\a^*}(2,0))\subset \pi_1(M_{i,\a^*}(2,0))$. Then there
is a lifting $K'_{\a^*}(2,0) \to M_i'$, by the covering lifting
theorem, and therefore we have a lifting $K'\to M_i''$ by lemma
\ref{branched lifting}. Because $im\{\pi_1(K)\to \pi_1(M_i'')\}
\subset im\{\pi_1(C)\to \pi_1(M_i'')\}$, we may then lift $K
\hookrightarrow N_i$. Thus, $M$ is a geometric limit of $N_i$. By
our construction, $N_i\simeq M$, and $C\subset K \subset M$ a
compact core, then $K\hookrightarrow N_i$ restricts to
$C\hookrightarrow N_i$ so that $C\simeq N_i$. Thus, $M$ is a super
strong limit of $N_i$.

Notice that $M_i'$ has indecomposable fundamental group, and
therefore is tame by Bonahon/Canary \cite{Bon, Ca89, Ca93}. Thus,
$M_i''$ is tame,  and therefore $N_i$ is tame by Cor. 3.2,
\cite{Ca94}.

\end{proof}

\section{Drilling and tameness} \label{sequence}
Recall that we have a PNC manifold $M$ and an algebraically
diskbusting link $\a^*\subset M$. We would like to have
$\pi_1(V_{\a^*}(2,0))\cong \pi_1(M_{\a^*}(2,0))$, but this might
not hold in general.

\begin{lemma} \label{tamereduction}
$M$ is a PNC manifold with $\pi_1(M)$ finitely generated. Let
$\a^*$ be an algebraically diskbusting geodesic link in $M$, $V$
an almost end-reduction of $M$ at $\a^*$. Then $V$ is tame, and
may be isotoped so that $V=\Int C$, where $C$ is a compact core of
$M$.
\end{lemma}
\begin{proof}
Take the orbifold cover $N\to M_{\a^*}(2,0)$ such that
$\pi_1(N)=\pi_1(V_{\a^*}(2,0))$. Then $V_{\a^*}(2,0)$ lifts
isometrically to $N$. Thus we may fill $\MN(\a^*)$ back in to get
a negatively curved manifold $N'$ which is a branched cover of $M$
and an isometric embedding $\i:V\hookrightarrow N'$ by lemma
\ref{branched lifting}.  Clearly $\i(V)$ is an almost
end-reduction at $\i(\MN(\a^*))\subset N'$. By construction,
$\pi_1 V_{\a^*}(2,0)=\pi_1 N$, so by theorem \ref{case1}, $N'$ is
tame. By Haken finiteness for $N'-\MN(\a^*)$, there exists $I$
such that $\d M_i$ is parallel to $\d M_j$ in $N'-\MN(\a^*)$ for
$i,j>I$, so that they also must be parallel in $V-\MN(\a^*)$. This
means that $V$ is tame and may be taken to be the interior of a
compact core $ C \subset M$, $V=\Int C$, by the uniqueness up to
isotopy  (theorem \ref{end reduction}).
\end{proof}

What we've just shown is that every algebraically diskbusting
geodesic link $\a^*$ is contained in a compact core $C$ of $M$.

\begin{theorem} \label{tame}
Let $M$ be PNC (without cusps), and $\pi_1(M)$ finitely generated.
Then $M$ is tame.
\end{theorem}
\begin{proof}
We may assume that if $C$ is a compact core for $M$, then $C$ is a
compression body. Thus all but one of the ends of $M$  will
correspond to incompressible surfaces in $\d C$, so by Canary's
theorem \cite{Ca93} these ends will be tame. Thus, there is a
single end $E$ corresponding to the compressible end of $\d C$
which we need to show is tame. There is a sequence of geodesics
$\b_i^*$ exiting $E$, by the argument of Theorem 3.10,
\cite{BrockBromberg02}. Now, if $\a\in \pi_1 M$ is an
algebraically diskbusting element, then $\a^*\cup \b_i^*$ is an
algebraically diskbusting link (possibly after perturbing the
metric on $M$ slightly). By lemma \ref{end reduction}, there
exists an end reduction $V_i$ at $\a^*\cup \b_i^*$. By lemma
\ref{Myers} $V_i$ is an almost end reduction at $\a^*\cup\b_i^*$.
Thus, there is a compact core $C$ of $M$ such that $\a^*\cup
\b_i^*\subset C$ by lemma \ref{tamereduction}. Since $\a$ is
algebraically diskbusting, $\d C$ is incompressible in $M-\a^*$.
But this implies that we may assume that $C\subset V$, where $V$
is the end-reduction of $M$ at $\a^*$, by the engulfing property.
Thus, $\b_i^*\subset V_{\a^*}(2,0)$ lifts to $V_{\a^*}(2,0)\subset
N$, where $N$ is the orbifold cover of $M_{\a^*}(2,0)$
corresponding to $\pi_1(V_{\a^*}(2,0))$, as in the proof of lemma
\ref{tamereduction}. Then the sequence of lifts of geodesics
$\b_i^*$ must exit the end $F$ of $N'$ corresponding to the
compressible component of $\d C$. Thus, the compressible end $F$
of $N'$ is geometrically infinite, and is therefore simply
degenerate by \cite{Ca93, Hou03}. By the PNC covering theorem
\ref{covering}, the end $F$ of $N$ must cover finite-to-one an end
$E$ of $M_{\a^*}(2,0)$. But the incompressible ends of $N'$ cover
one-to-one the incompressible ends of $M$, which means that $F$
must cover the compressible end $E$ finite-to-one, and therefore
$E$ is tame by \ref{finitecovering}.
\end{proof}

\section{Isoperimetric inequality}

Let $\widetilde{\CH(M)}\subset \widetilde{M}$ be the universal
cover of $\CH(M)$, where $(M,g)$ is a PNC manifold. Let $s\subset
\d \widetilde{\CH(M)}$ be a simple closed curve which is
contractible in $\d \widetilde{\CH(M)}$. Then $s=\d D$, $D\subset
\widetilde{\CH(M)}$.

\begin{lemma} \label{isoperimetric}
There exists a constant $r$ depending only on the pinching
constants, such that $D\subset \MN_{r+diam(s)/2}(s)$.
\end{lemma}
\begin{proof}
Consider $\CH(s)\subset \widetilde{\CH(M)}$. We claim that
$D\subset \d \CH(s)$. Since $\CH(s)$ is a ball, and $s\subset \d
\CH(s)$ (since $s\subset \d \widetilde{\CH(M)}$), then $s$
separates $\d \CH(s)$ into two disks $\d \CH(s)= E_1\cup E_2$, $\d
E_i=s$. Let $E_1$ be the disk which is closer to $D$. Then
$E_1\cup D$ is the frontier of a closed contractible subset $B
\subset \widetilde{\CH(M)}$. Then $U=\widetilde{\CH(M)}-B\cup E_1$
is a closed convex set. If not, then there is a pair of points
$x,y\in U$ and a geodesic $[x,y]$ connecting $x$ and $y$, such
that there is a subinterval $[x_1,y_1]\subset \widetilde{M}-\Int
U$ (since $E_1$ is separating in $\widetilde{\CH(M)}$). Then
$x_1,y_1\subset \d U$, so they must lie in $E_1$, otherwise we
would violate the convexity of $\widetilde{\CH(M)}$. But
$E_1\subset \CH(s)$ which is convex, so $[x_1,y_1]\subset
\CH(s)\subset U$, a contradiction. Thus,
$\widetilde{\CH(M)}\subset U$, so $D=E_1$, as claimed.

For a set $Q$, let $join(Q)$ be the union of all geodesics
connecting pairs of points in $Q$. Clearly $join(Q)\subset
\MN_{diam(Q)/2}(Q)$. Now, by the arguments on pp. 241-243 of
\cite{Bowditch95}, $\CH(s)\subset \MN_r(join(s))$, for some $r$
which depends only on the pinching constants of $(M,g)$. Thus, we
have $D\subset \CH(s)\subset \MN_{r+diam(s)/2}(s)$.

\end{proof}

\section{Interpolation of simplicial ruled surfaces} \label{interpolate}

The following arguments are modifications of work of Brock in
\cite{Bro01}, as a way of generalizing the simplicial ruled
surface interpolation techniques of Thurston \cite{Th} and Canary
\cite{Ca94}. We were unable to generalize Canary's method directly
from the hyperbolic case to the PNC case. Let $\S$ be a closed
surface of negative Euler characteristic. A {\it pants
decomposition} $P$ of $\S$ is a maximal collection of isotopy
classes of disjoint non-parallel essential simple closed curves on
$\S$. Two distinct pants decompositions $P$ and $P'$ are related
by an {\it elementary move} if $P'$ can be obtained from $P$ by
replacing a curve $\a\in P$ by a curve $\b\in P'$ intersecting
$\a$ minimally (and disjoint from the other curves of $P$, {\it
i.e.} so $\a\cap P'=\b\cap P = \a\cap \b$, and $ |\a\cap \b|=1$ if
$\a$ is non-separating in the component of $\S-(P-\a)$ containing
$\a$, and $|\a\cap\b|=2$ if $\a$ is separating in the component of
$\S-(P-\a)$ containing $\a$). The {\it pants graph} $P(\S)$ is the
graph with pants decompositions of $\S$ as vertices, and edges
joining pants decompositions which differ by an elementary move.

Given a pants decomposition $P$ of $\S$, we may extend the
collection of curves to a {\it standard triangulation}. These
triangulations are obtained by adding a pair of vertices to each
curve of $P$, then extending to a triangulation of each pants
region by adding a maximal collection of non-parallel edges (we
allow any extension, unlike Brock).

Given a pants decomposition $P$ of $\S$, a standard triangulation
$T$ extending $P$, and a $\pi_1$-injective map $f:\S\to N$, where
$N$  is a PNC orbifold, we may homotope the map $f:\S\to N$ to a
simplicial ruled surface. First, we homotope the loops of $P$ to
be geodesic, then we homotope the other edges of $T$ to be
geodesic, then we simplicially interpolate the triangles of $T$.
This gives a simplicial ruled surface $f^*:\S\to N$, by prop.
3.2.7 \cite{Ca89}.

Given a fixed pants decomposition $P$, two standard triangulations
$T_1,T_2$ associated to $P$, and simplicial ruled representatives
$f_1^*:\S\to N$, $f_2^*:\S\to N$  in a homotopy class $f:\S\to N$,
we may interpolate continuously between these simplicial ruled
representatives. This is done as in Claim 3 of the proof of
theorem 2.3, \cite{BCW04}. The moves consist of changing the
ruling on a triangle or doing a flip move between triangles which
are adjacent in $\S-P$. We also need a move which moves a vertex
along a geodesic representative of a curve of $P$, so that we can
get the vertices to line up. All of these moves preserve the
simplicial ruled conditions.

We also need moves to get from $P$ to $P'$, where $P$ and $P'$ are
adjacent pants decompositions. We may choose simplicial ruled
representatives $f_1^*:\S\to N$, $f_2^*:\S\to N$, where $f_1^*$
has a standard triangulation with respect to $P$, and $f_2^*$ has
a standard triangulation with respect to $P'$.

Claim: We may continuously move $f_1^*$ and $f_2^*$ through
simplicial ruled surfaces to $f_1'$ and $f_2'$, so that $f_1'$ is
homotopic to $f_2'$ in $\MN_R(f_1'(\S)\cup f_2'(\S))$, where $R$
only depends on the pinching constants.

\begin{figure}[htb] \label{pants}
    \begin{center}
    \epsfbox{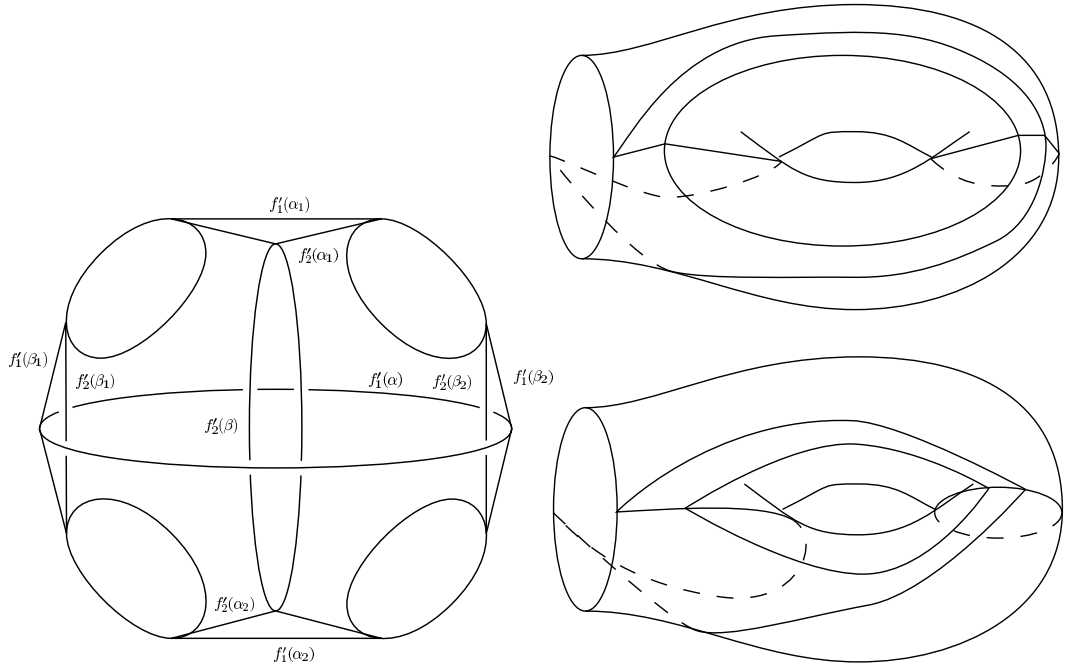}
    \caption{\label{triangulation} Interpolation of standard triangulations}
    \end{center}
\end{figure}

We need only show how the triangulations behave on the 4-punctured
sphere and punctured torus regions; we may assume that the
standard triangulations agree on the rest of the pair of surfaces,
by using the previously described moves to make the standard
triangulations line up on this subsurface. For the 4-punctured
sphere, insert arcs $\a_1,\a_2$ disjoint from  $\a$ and
$\b_1,\b_2$ disjoint from $\b$ into both 4-punctured spheres in
$\S$ so that the corresponding edges in $f_1'(\S)$ and $f_2'(\S)$
are homotopic in $M$ rel endpoints, then complete this to a
triangulation. $\a\cup \a_1\cup \a_2$ will be geodesic in $f_1'$,
so $\b_1\cup\b_2$ will be composed of two geodesics (since these
meet $\a$ once). In the other surface $f_2'$, $\b_1\cup\b_2$ will
be geodesic. Together, $f_1'(\b_i)\cup f_2'(\b_i)$ will form a
geodesic triangle, and similarly for $f_1'(\a_i)\cup f_2'(\a_i)$,
$i=1,2$ (see figure \ref{pants}).

There are four such triangles, which we may assume are simplicial
ruled. The triangles will lie in a uniform $\delta$ neighborhood
of the edges, by $\delta$-hyperbolicity (where $\delta$ only
depends on the pinching constants). If we take four triangles, and
one half of each 4-punctured sphere, we get a map $f:S^2\to M$ of
a sphere of bounded area (depending only on the pinching
constants). Let $\tilde{f}:S^2\to \tilde{M}$ denote a lift of $f$
to $\tilde{M}$.

\begin{lemma}
Let $\tilde{M}$ be a PNC simply connected manifold.  Let
$\tilde{f}:S^2\to \tilde{M}$, then $\tilde{f}$ is homotopically
trivial in $\MN_R(\tilde{f}(S^2))$, where $R$ depends only on the
area of $\tilde{f}$.
\end{lemma}
\begin{proof}
If we add the bounded complementary regions of $\tilde{f}(S^2)$,
and take a small regular neighborhood to get a submanifold
$U\subset \tilde{M}$, so that $\tilde{M}-U$ has no bounded
components, then $\pi_2(U)=0$ by standard arguments, so
$\tilde{f}$ is homotopically trivial in $U$. Thus, if $\tilde{f}$
is not homotopically trivial in $\MN_R(\tilde{f}(S^2))$, this
means that $\tilde{M}-\MN_R(\tilde{f}(S^2))$ must have a bounded
component $V$ which is contained in a bounded component $V'$ in
the complement of $\tilde{f}(S^2)$. Choose $x\in V$, then
$d(x,\tilde{f}(S^2))> R$. But $\d V'\subset \tilde{f}(S^2)$ has
$\area(\d V') < \area(\tilde{f}(S^2))$. Since $\vol(V')>
\vol(B_R(x))$, by the isoperimetric inequality for $\tilde{M}$
({\it e.g.} \cite{Kleiner92}), this implies that
$\area(\tilde{f}(S^2))>\area(\d V')> \area(\d B_R(x))$. Thus,
choosing $R$ large enough that $\area(\d B_R(x))>
\area(\tilde{f}(S^2))$, we get a contradiction.
\end{proof}

By this lemma, $f$ is homotopically trivial in $\MN_R(f(S^2))$, by
projecting $\MN_R(\tilde{f}(S^2))$ to $M$. Thus $f_1'$ and $ f_2'$
will be homotopic in $\MN_R(f_1'(\S)\cup f_2'(\S))$. A similar
construction works for elementary moves involving punctured tori.
Since in $\MN_R(f_1'(\S)\cup f_2'(\S))$, there is a homotopy
between the two surfaces, every point in the support of this
homotopy will be within a distance $R$ of a simplicial ruled
surface. A theorem of Hatcher, Lochak, and Schneps \cite{HLS00}
implies that the pants complex is connected. This implies that we
may interpolate between any two simplicial ruled surfaces
realizing pants decompositions by surfaces which are bounded
distance from simplicial ruled surfaces.

\section{Tameness of limits} \label{limits}

The following is a generalization of the main theorem of
\cite{CM96}.
\begin{lemma}
Suppose $N_i$ are tame PNC manifolds (with no cusps) which super
strongly converge to $M$. Then $M$ is tame.
\end{lemma}

\begin{proof}
We may restrict to the case that $M$ has a compression body core
$C$, since this is the case that we use in the application to
tameness. We may assume that the compressible end of $M$ is not
geometrically finite, otherwise we are done.

The first step of the argument is to show that $\d \CH(N_i)\to
\infty$, in the appropriate sense (Dick Canary suggested this line
of argument, which gives a simplification over our previous
approach). For simplicity, we will assume that $M$ has only one
end. To modify the argument appropriately, replace $\CH(N_i)$ with
the convex submanifold $\CH(N_i)\cup E$, where $E$ consists of the
ends of $M-\CH(N_i)$ corresponding to the incompressible
components of $\d C$.

Claim: for all compact connected $K\subset M$, there exists $I$
such that for $i>I$, $\eta_{i,K}(K)\subset \CH(N_i)$ (where the
maps $\eta_{i,K}$ come from the definition \ref{strong limit} of
super strong convergence ).

Fix a compact core $C$ for $M$.  Suppose not, then there is a
connected $K\subset M$ and $J\subset \NN$ with $|J|=\infty$, such
that $\eta_{j,K}(K)\nsubseteq \CH(N_j)$, for all $j\in J$. We may
choose an algebraically diskbusting geodesic $\a^*\in M$, and
assume $C\cup \MN(\a^*)\subset K$ by enlarging $K$. Then there is
an $I$ such that for $i\geq I$, we have isometric maps
$\eta_{i,\a^*}:\a^* \to N_i$.

 Then there is a $\delta$ depending on $\e$, the
Margulis constant, such that $diam_{\e}(\d \CH(N_j)) \leq A/\e +
\delta$, for $j\in J$.

The region $\MN^{\e}_R(\eta_{j,K}(\a^*))$ is a finite union of
tubular neighborhoods of Margulis tubes and geodesics. Thus, we
may perturb $\d \CH(N_j)$ slightly to meet $\d
\MN^{\e}_R(\eta_{j,K}(\a^*))$ transversely for generic $R$ in a
collection of simple curves. Consider the interval $[R_1,R_2]$
such that $L(R)=length(\d \MN^{\e}_R(\eta_{j,K}(\a^*))\cap \d
\CH(N_j))
> \e$, for generic $R\in [R_1,R_2]$, and some component of $\d
\MN^{\e}_R(\eta_{j,K}(\a^*))\cap \d \CH(N_j)$ is homotopically
non-trivial in $N_j$.

The components of $\d \MN^{\e}_{R_1}(\eta_{j,K}(\a^*))\cap \d
\CH(N_j)$ which are trivial in $N_j$ but not in $\d \CH(N_j)$ must
have length $> f(R_1)$, where $f$ is a monotonic function
depending only on the pinching constants. This follows since any
disk they bound must intersect $\a^*$, so by the isoperimetric
inequality for a PNC manifold, the length must be greater than
that of the boundary of a disk of radius $R_1$ (the argument here
is the same as in Lemma 4.1 \cite{CM96}). So by declaring that
$R_1> f^{-1}(\e)$, we may assume that all components of $\d
\MN^{\e}_{R}(\eta_{j,K}(\a^*))\cap \d \CH(N_j)$ of length $<\e$
which are trivial in $N_j$ are actually trivial in $\d\CH(N_j)$
for $R\in [R_1,R_2]$.

The reason the interval $[R_1,R_2]$ exists is that if for some
$R$, $L(R) \leq \e$, then either a component of  $\d
\MN^{\e}_{R}(\eta_{j,K}(\a^*))\cap \d \CH(N_j)$  is homotopically
non-trivial in $N_j$ and has length $<\e$, in which case it lies
in $M_{thin(\e)}$, and therefore cannot lie in $\d
\MN^{\e}_R(\eta_{j,K}(\a^*))$ for generic $R$, or every component
has length $<\e$ and is homotopically trivial in $\CH(N_j)$, in
which case either every loop before or every loop after is
homotopically trivial (this follows from an outermost loop
argument). Then $\e (R_2-R_1) \leq \int_{R_1}^{R_2} L(r) dr \leq
A$ (by the coarea formula). Thus, $R_2-R_1\leq A/\e$.

 We know that
$R_1$ cannot be too large. Every loop in $\d
\MN^{\e}_{R_1}(\eta_{j,K}(\a^*))\cap \d \CH(N_j)$ is trivial in
$\d \CH(N_j)$, and therefore bounds a disk of diameter $<r+\e$ by
lemma \ref{isoperimetric}. Therefore $\d \CH(N_j)$ lies outside of
$\MN^{\e}_{R_1-r-\e}(\eta_{j,K}(\a^*))$, which would contradict
the fact that $\d \CH(N_j)$ meets $K$ if $R_1$ were too large.
Similarly, $R_2$ cannot be too large, since $R_2\leq R_1+A/\e$. We
also know that all of the loops of $\d
\MN^{\e}_{R_2}(\eta_{j,K}(\a^*))\cap \d \CH(N_j)$ bound disks in
$\d \CH(N_j)$ of diameter $< r+ \e$. Thus, $\d \CH(N_j)\subset
\MN^{\e}_{R(K)}(\eta_{j,K}(\a^*))$, where $R(K)$ depends on our
set $K$.

Thus, if $K\cap \d \CH(N_j)\neq \emptyset$ for all $j\in J$, then
$\d \CH(N_j)\subset \MN^{\e}_{R(K)}(\eta_{j,K}(\a^*))$. Now, we
may take $Q=\MN^{\e}_{R(K)}(\a^*)\subset M$, and for $j$ large
enough, this maps isometrically to $N_j$ by $\eta_{j,Q}$, and
therefore $\eta_{j,Q}(Q)=\MN^{\e}_{R(K)}(\eta_{j,Q}(\a^*))$ in
$N_j$. Therefore, $\CH(N_j)\subset \eta_{j,Q}(Q)$. Pushing back to
$M$, we see that $M$ is geometrically finite, contradicting our
assumption. This finishes the proof of the claim.

The argument of Canary-Minsky \cite{CM96} proceeds by finding
simplicial hyperbolic surfaces which are arbitrarily close to $\d
\CH(N_i)$. We will work a little more coarsely, by finding a
simplicial ruled surface which has points mod $\e$-close to $\d
\CH(N_i)$ (this line of argument was suggested by Minsky).

Take a systole $g\subset \d_0 \CH(N_i)=S$ (the compressible
component of $\d\CH(N_i)$). By theorems \ref{boundary} and
\ref{systole}, there is an $l$ depending only on the pinching
constants and $\chi(M)$ such that $length(g)\leq l$. Then we may
homotope $g$ to a geodesic $g^*$ in $\MN_{R}^{\e}(\d_0 \CH(N_i))$,
for some $R$ depending only on $\chi(M)$ and the pinching
constants. If $R$ is large enough, then $g$ will be homotopic to
$g^*$ in the complement of $\MN(\eta_{i,\a^*}(\a^*))$.  We may
complete $g$ to a pants decomposition $P$, and a standard
triangulation of $S$, to get an incompressible simplicial ruled
surface $f:S\to N_{i,\a^*}(2,0)$ realizing $P$ as geodesics. We
may choose a minimal path from $P$ to a pants decomposition $P'$
so that $P'$ has a compressible curve. Then we may use the moves
described in section \ref{interpolate} to interpolate between
simplicial ruled surface representatives in $N_{i,\a^*}(2,0)$
(where $\a^*$ is a fixed algebraically diskbusting knot in $M$,
which we may map isometrically into $N_i$ for $i$ large enough)
with standard triangulations associated to $P$ and $P'$. The
geodesics representative of $P'$ must meet $\MN(\a^*)$, since
otherwise the compressible curve of $P'$ would have a geodesic
representative in $N_i$.

One surface in this homotopy will meet $\MN^{\e}_R(\a^*)$, in
which case there is a simplicial ruled surface within distance
$R$. We can push these simplicial ruled surfaces to
$M_{\a^*}(2,0)$, to get longer and longer sequences of simplicial
ruled surfaces exiting the (distinguished) end of $M_{\a^*}(2,0)$
(if the compressible end of $N_i$ is simply degenerate, then we
get these surfaces from the filling theorem \ref{filling}). We may
homotope these surfaces in the complement of $\a^*$ by section
\ref{interpolate} to simplicial ruled surfaces lying in a compact
subset $K$ of $M_{\a^*}(2,0)$. By Souto's finiteness theorem
(prop. 2, section 2.1, \cite{Souto02}), these surfaces are
eventually homotopic (in $M_{\a^*}(2,0)$), and so we see that we
have a homotopy of a surface properly exiting the compressible end
of $M_{\a^*}(2,0)$, such that every point in the image of this
homotopy is within distance $R$ of a simplicial ruled surface.
Then we may finish off using the argument of the main theorem
section 8 \cite{CM96} or theorem 1, section 3 of \cite{Souto02} to
show that $M$ is tame.
\end{proof}

\section{PNC covering theorem}

Canary proved a covering theorem for simply degenerate ends of
hyperbolic manifolds \cite{Ca94}. The main element of Canary's
argument that needs to be generalized to the PNC category is the
interpolation of pleated surfaces, which is described in section
\ref{interpolate}. In fact, for our application, we only need to
apply the covering theorem for incompressible ends (this was
proven by Thurston \cite{Th} in the hyperbolic category).

\begin{theorem} \label{filling} (the Filling theorem)
Let $N$ be a PNC tame 3-orbifold without cusps, and $E$ a simply
degenerate incompressible (manifold) end of $N$. There exists a
constant $R$ such that $E$ has a neighborhood $U$ homeomorphic to
$S\times [0,\infty)$ such that every point in some subneighborhood
$\hat{U}=S\times [k,\infty)$ is distance $\leq R$ from a
simplicial ruled surface $f:S\to \hat{U}$.
\end{theorem}
\begin{proof}(Sketch)
 The idea is straightforward here in the
incompressible case. Because the end is simply degenerate, we may
realize one simplicial ruled surface using a standard
triangulation associated to a pants decomposition. Since the end
is simply degenerate, we may find simple closed curves on $S$ with
geodesic representatives exiting the end. We extend these curves
to pants decompositions, and standard triangulations, and then use
section \ref{interpolate} to homotope between these and the first
surface with surfaces of bounded mod $\e$ diameter, so that every
point in the support of the homotopy is bounded distance from a
simplicial ruled surface. Then the filling theorem follows from
standard topological arguments and the fact that the end has
$\infty$ mod $\e$ diameter by lemma \ref{diameter}.
\end{proof}

\begin{theorem} (the PNC Covering theorem) \label{covering}
Let $p:M\to N$ be an orbifold cover, where $N$ is PNC and $M$ is
tame. If $E$ is a simply degenerate incompressible manifold end of
$M$ (\it{i.e.} $E\cong S\times [0,\infty)$ for some surface $S$),
then either
\begin{enumerate}
\item
$E$ has a neighborhood $U$ such that $p$ is finite-to-one on $U$,
or
\item
$N$ has finite volume and has a finite cover $N'$ which fibers
over the circle, such that if $N_S$ denotes the cover of $N'$
associated to the fiber subgroup, then $M$ is finitely covered by
$N_S$. Moreover, if $M\neq N_S$, then $M$ is homeomorphic to the
interior of a twisted $I$-bundle which is double covered by $N_S$.
\end{enumerate}
\end{theorem}
 \begin{proof}(sketch) This theorem follows Canary's argument
closely. If the map of $p_{|E}$ is not finite-to-one, then by the
filling theorem \ref{filling}, there is a point $x\in N$, and
infinitely many simplicial ruled surfaces $p\circ f_i:S\to N$
meeting $B_R(x)$, where $f_i:S\to E$ are simplicial ruled surfaces
exiting $E$. By the finiteness theorem (see Prop. 2 \cite{Souto02}
for a PNC version), infinitely many of these surfaces must be
homotopic. By theorem 2.5 \cite{CM96}, we may find embedded
surfaces $S_i\subset \MN(f_i(S))\subset E$, so that
$(S_i\hookrightarrow N)\simeq f_i$. But if $p\circ S_i$ and
$p\circ S_j$ are homotopic in a bounded neighborhood of $x\in N$,
for $j\gg i$, then this means we have an embedded lift
$f_i':S_i\to E$ so that $f_i'(S_i)$ is near to $S_j$, and
$p(S_i)=p( f_i'(S_i))$. But then the product region between $S_i$
and $f_i'(S_i)$ gives a mapping torus $S_i\times [0,1]/ \{(s,0)
\sim (f_i'(s),1)\}$, which is a finite sheeted fibered cover of
$N$.
\end{proof}

We need a slightly stronger result than that provided by Canary's
theorem:

\begin{lemma} (finite covering theorem) \label{finitecovering}
Let $M\to N$ be an (orbifold) cover, and let $E\subset M$,
$F\subset N$ be ends, such that $E\to F$ is finite-to-one. Suppose
that $E$ is a tame manifold end, so there is a neighborhood
$E'\subset E$ so that $E'\cong S\times [0,\infty)$. Then $F$ is
tame as well, so there is a 2-orbifold $\MO$ and a neighborhood
$F'\subset F$ such that $F'\cong \MO\times [0,\infty)$.
\end{lemma}
\begin{proof}
By replacing $E$ with a finite sheeted cover of $E$, we may assume
that the cover $E\to F$ is regular. We may also assume that there
is a product neighborhood $E_0$ such that $E\subset E_0$. Let
$E_2\subset E=E_1$ be a product neighborhood, and let $F_3\subset
F_1=F$ be a closed end neighborhood so that $\d F_3$ is a
suborbifold, and $F_3$ is covered by $E_3\subset E_2$. Then $\d
E_3$ is a surface covering $\d F_3$. Suppose that $\d E_3$ is
compressible in $E_2$. By the equivariant Dehn's lemma, we may
find a set of compressing disks for $\d E_3$ in $E_1$ which are
invariant with respect to the covering translations of $E_1\to
F_1$, and thus passes to an orbifold compression of $\d F_3$. Do a
maximal set of such compressions to get a closed neighborhood
$F_3'$ covered by $E_3'\subset E_1\subset E_0$, so that  $\d E_3'$
is incompressible in $E_0$. Since $E_1$ is a product, $\d E_3'$
must be a surface parallel to $\d E_1$, since the only
incompressible surface in a product is boundary parallel.
Continuing in this manner, we get a sequence of equivariant
product neighborhoods $E_4, E_5, ...$ exiting the end. In
$E_i-\Int E_{i+1}$, the action is standard, so the covering action
on $E_i$ must be standard, and we see that $F$ is also a product.

\end{proof}

\appendix
\section{Boundary of the convex core}
The goal of this appendix is to give details of Kleiner's argument
on pp. 42-43 of \cite{Kleiner92}, which we use to bound the area
of the boundary of the convex core.  Let $M^3$ be a tame PNC
manifold, without cusps, and pinching constants $a< K(P) <b<0,
\forall P \in G_2(TM)$. Assume that $dim(\CH(M))=3$.
\begin{theorem} \label{boundary}
There is a constant $A<0$ which depends only on the pinching
constants such that $\area(\d \CH(M))<A \chi(\d \CH(M))$.
\end{theorem}

\begin{proof}
By geometric tameness, $\d \CH(M)$ will be a compact surface
\cite{Ca89, Hou03}. The philosophy of the argument is
straightforward. If $\d \CH(M)$ were smooth, then at each point of
$\d \CH(M)$, the Gauss-Kronecker curvature would be zero,
otherwise we could push in slightly and get a smaller convex
submanifold. Then the Gauss-Bonnet theorem would bound the area of
$\d \CH(M)$. Since $\d \CH(M)$ is not necessarily smooth, we
perform this argument for nearby surfaces which are smooth enough,
then take a limit, showing that the average Gauss-Kronecker
curvature approaches zero.

 We follow the argument of
Kleiner \cite{Kleiner92}, including some more details. Kleiner was
taking the convex hull of a compact set, but his argument is
local, so it generalizes to the convex hull of the limit set of
$\pi_1 M$.

Let $\delta: M-\CH(M)\to (0,\infty)$ be the function such that
$\delta(x)=d(x,\CH(M))$. Since balls in a negatively curved
manifold are strictly convex ({\it e.g.} see section 7.6
\cite{Cha93}), for any $x\in M-\CH(M)$, $B_{\delta(x)}(x)\cap
\CH(M)=\xi_0(x)$, for a unique point $\xi_0(x)\in \d \CH(M)$ (see
1.6 \cite{BGS85}, working equivariantly with
$\tilde{M}-\widetilde{\CH(M)}$). Let $E_s$ be the points in $M$
distance $\leq s$ away from $\CH(M)$, so
$E_s=\delta^{-1}(0,s]\cup\CH(M)$. Let $C_s=\d E_s=\delta^{-1}(s)$.
We see that $E_s$ is convex, since $\delta$ is a convex function
by 1.6(iii) \cite{BGS85}. There is also a nearest point retraction
$\xi_s:M-E_s\to C_s$. These functions are 1-Lipschitz, as shown in
1.6(i) \cite{BGS85}.

\begin{lemma}
The vector field $\nabla \delta$ is locally Lipschitz on $M-E_s$,
for any $s>0$.
\end{lemma}

\begin{proof}
This argument follows Federer (pp. 433-440, 4.7 \& 4.8
\cite{Fed76}), which proves a corresponding result for convex
hulls in Euclidean spaces. We can work in $\tilde{M}$ since we
need only a local result, but for simplicity we will abuse
notation and use the same notation as for $M$. The function
$\delta$ is clearly 1-Lipschitz, by the triangle inequality. In
fact, $\delta$ is differentiable, which follows essentially from
the continuity of $\xi_0(x)$ (this follows from the argument
below, and see Lemma 4.7 \cite{Fed76}). $\nabla \delta$ is a unit
vector field on $M-\CH(M)$, orthogonal to the foliation $\{C_s\}$.
$\nabla\delta$ is $C^{\infty}$ in the direction $\nabla\delta$,
since the flow lines of $\nabla \delta$ are length-minimizing
geodesics. Intuitively, since the surfaces $C_s$ are supported by
balls of the form $B_{s}(\xi_0(x)), x \in C_s$, then
$\nabla\delta_{|C_s}$ should have derivative bounded by the
curvature of $\d B_s(\xi_0(x))$.

To make this intuition precise, we compare the vector field
$\nabla\delta $ with radial vector fields. Fix $s_0>0$, then for
any point $x\in M-E_{s_0}$, $V_x=\nabla d(\cdot,\xi_0(x))$ is a
$C^\infty$ (unit) vector field with $\nabla V_x$ bounded
independent of $x$ (but dependent on $s_0$). This follows since
$\nabla V_{x| \d B_r(\xi_0(x))}$ is the shape operator for $\d
B_r(\xi_0(x))$, which has bounded norm (dependent on $s_0$)
following from theorem 3.14 \cite{Cha93}. We need to show that
$\nabla\delta $ is locally Lipschitz. Fix some bounded open
neighborhood $U\subset M-E_{s_0}$ such that $diam(U)\leq \e$, and
a smooth chart $U\hookrightarrow\RR^3$. The Euclidean metric
induced on $U$ is bi-Lipschitz to the metric on $U\subset M$. Then
there is a constant $C$, such that for any $x,y \in U$,
$|V_x(x)-V_x(y)|< C |x-y|$. Now, consider distances
$a_1=d(y,\xi_0(x)), a_2=d(y,\xi_0(y)), a_3=d(\xi_0(x),\xi_0(y))$.
Then we have $a_3\leq d(x,y)\leq \e$. By the Toponogov comparison
theorem (see 1.5 \cite{BGS85}), we have $\langle V_x(y),
\nabla\delta(y)\rangle \geq
\frac{a_1^2+a_2^2-a_3^2}{2a_1a_2}=\frac12(a_1/a_2+a_2/a_1 -
a_3^2/a_1a_2) \geq 1-\e^2/2s_0^2$. Thus,
$|V_x(y)-\nabla\delta(y)|\leq \e/s_0$. Since
$\nabla\delta(x)=V_x(x)$, then we obtain
$|\nabla\delta(x)-\nabla\delta(y)|\leq
|V_x(x)-V_x(y)|+|V_x(y)-\nabla\delta(y)| =
O(|x-y|)+O(d(x,y))=O(|x-y|)$. Thus, $\nabla\delta$ is locally
Lipschitz. Moreover, if $\nabla\delta$ is $C^2$ at $x$, then
$\nabla\delta$ is $C^2$ along the ray
$\overrightarrow{\xi_0(x)x}$, since solutions to the Jacobi
equation will be $C^2$. By Rademacher's theorem, $\nabla\delta$ is
$C^2$ a.e. Also, $\nabla\delta_{|C_s}$ is Lipschitz, and therefore
$C^2$ a.e.

\end{proof}

Let $GK_{C_s}$ denote the Gauss-Kronecker curvature of $C_s$, the
product of the principal curvatures, that is the determinant of
$A^{\nu}_{C_s}$, the shape operator. Since $\nabla\delta_{|C_s}$
is $C^2 a.e.$, this implies that $GK_{C_s}$ is a well-defined
element of $L^{\infty}(C_s)$.

Claim: $\int_{C_s} GK_{C_s} da \to 0$ as $s\to 0^+$.

  For $0 \leq
s \leq s_0$, let $r_{s_0s}:C_{s_0}\to C_s$ be the Lipschitz
nearest point projection (so $r_{s_0s}=\xi_{s|C_{s_0}}$). Assume
that $\nabla\delta$ is $C^2$ at $p\in C_{s_0}$, and
$\k=\inf\{K(P)\}$ where $P$ runs over all 2-planes. We'll define
$c_{\k}(t)=\cosh(\sqrt{-\k t})$, $s_{\k}(t)=(-\k)^{-\frac12}
\sinh(\sqrt{-\k t})$.  To prove the claim, we need to show that

\begin{enumerate}
\item $|GK_{C_s}(r_{s_0s}(p)) Jac(r_{s_0s}) (p)|\leq F(s_0,\k)$,
where $F$ is some function,

\item
$\underset{s\to 0}{\lim} ( GK_{C_s}(r_{s_0s}(p)) Jac(r_{s_0s})(p)
= 0$.
\end{enumerate}

Then $$\int_{C_s} GK_{C_s} d\area_{C_s} = \int_{C_{s_0}}
(GK_{C_s}\circ r_{s_0s}) Jac(r_{s_0s}) d\area_{C_{s_0}} \to 0$$ as
$s\to 0$. Let $\nu=-\nabla\delta$. Pick $p\in C_{s_0}$ at which
$\nu$ is differentiable. Let $\g:[0,s_0]\to M$ be the geodesic
segment $\g(t)=\exp t\nu(p)$, and for every $e\in T_p C_{s_0}$,
let $Y$ be the Jacobi field along $\g$ given by
$Y(t)=(\exp\circ(t\cdot \nu))_* e$. Thus, $Y(0)=e$.

For $s\in (0,s_0]$ consider maps $W_{s_0s}:T_pC_{s_0}\to
T_{r_{s_0s}(p)}C_s$ given by $e\to \nabla_{\dot{\g}(s_0-s)}Y$. We
claim that the maps $W_{s_0s}$ are bounded above uniformly in
terms of $s_0$ and the geometry of $M$, while the lower bound on
$W_{s_0s}$, {\it i.e.} $\inf\{|W_{s_0s} e| | e\in T_pC_{s_0},
|e|=1\}$, goes to zero as $s\to 0$. To see the former, note that
the second fundamental form of $C_{s_0}$ is bounded above
uniformly in terms of $s_0$ and the geometry of $M^3$ since
$C_{s_0}$ is convex and supported from the inside by a ball of
radius $s_0$; consequently the maps $W_{s_0s}$ are bounded above
uniformly because they are obtained by solving the Jacobi equation
with initial conditions determined by the second fundamental form
of $C_{s_0}$.

Assume that $|Y(0)|=1$. Then
$$\frac{d}{dt} \langle \nabla_t Y, \nabla_t Y \rangle = 2\langle
Y', Y''\rangle =$$
$$ -2 \langle Y', R(\g',Y)\g'\rangle = -2
|Y'||Y| R(\g',Y/|Y|,\g',Y'/|Y'|\rangle\leq -2\k |Y'||Y|.$$

Assuming $\l$ is the minimum eigenvalue of $A^{\nu}$, we have
$$ \frac{d}{dt} |Y'| = \frac{1}{2|Y'|} \frac{d}{dt} \langle
Y',Y'\rangle \leq -\k |Y| \leq -\k (c_{\k} -\l s_{\k})(t) |Y(0)|
\leq -\k c_{\k}(t),$$ since $\l\geq 0$ and $|Y(0)|=1$. This
follows from theorem 7.4, \cite{Cha93}. Integrating, we obtain
$$|Y'|(t) \leq \int_{0}^{s_0-t} -\k c_{\k}(u) du +|Y'(0)| = -\k
s_{\k}(s_0-t) + |Y'(0)|.$$

Thus, we need to bound $|Y'(0)|$. We have the formula
$Y'(0)=-A^{\nu}(Y(0))$ (p. 321 \cite{Cha93}). Intuitively,
$|A^{\nu}|$ is bounded, since the surface $C_{s_0}$ is supported
by the ball $B_{s_0}(\xi_0(p))$, which has principle curvatures
controlled purely by $\k$ and $s_0$. By theorem 3.14,
\cite{Cha93}, if $\tr A^{\nu}\geq (n-1)\l$, (in our case,
$n-1=2$), then the map $\exp$ has a focal point at distance the
smallest zero of the function $(c_{\k}-\l s_{k})(t)$. But $\exp$
has no focal point at distance $s_0$, so this means that $\tr
A^{\nu} \leq 2 c_{\k}(s_0)/s_{\k}(s_0)$. Since $C_{s_0}$ is
convex, the principle curvatures $\l_1,\l_2\geq 0$, so $\tr
A^{\nu}=\l_1+\l_2\geq \l_i$. Thus, $\l_i\leq 2
c_{\k}(s_0)/s_{\k}(s_0)$. Since we are assuming $|Y(0)|=1$, and
$|A^{\nu}|=\max\{\l_1,\l_2\}$, we have $|Y'(0)|\leq 2
c_{\k}(s_0)/s_{\k}(s_0)$, which clearly only depends on $\k$ and
$s_0$.

We also need to show that the lower bound $\inf \{|W_{s_0s}e|\
\mid e\in T_p C_{s_0}, |e|=1\}$ goes to zero as $s\to 0$. Note
that we have the factorization $W_{s_0s}= A^{\nu}_{C_s}\circ
r_{s_0s}^*$, where $A^{\nu}_{C_s}:T_{\g(s_0-s)}C_s\to
T_{\g(s_0-s)}C_s$ is the Weingarten map for the inward normal to
$C_s$ (given by $A^{\nu} u=-(\nabla_u \nu)^T$, where $u\in T_p
C_s)$), and the fact that the lower bound on $A^{\nu}$ goes to
zero with $s$.

 We prove this by contradiction. Suppose $A^{\nu}$
does not approach 0. Then there exists a sequence $s_i \to 0$ such
that $\inf\{|A^{\nu}_{s_i}e| | |e|=1,       e\in T_{r_{s_0s}(p)}
C_s\} \geq \e$. This is also equal to $\inf\{\l_1,\l_2\}$, where
$\l_i$ are the eigenvalues for $A^{\nu}_s$. In $T_{\g(s_0-s_i)}M$,
take a little spherical cap $C_i$ with radius of curvature $\e/4$
and tangent to $C_{s_i}$. At $\g(s_0-s_i)$, $\l_i(C_i)=\e/4$,
since the connections agree to first order at the origin. Since
$\l_i(C_{s_i})\geq \e$, then we have $\l_i(C_{s_i})\geq \e/2$ in a
small neighborhood of $\g(s-s_i)$ in $C_{s_i}$. Assume that the
spherical cap $C_i$ is small enough so that $0\leq \l_j(\exp
C_i)\leq \e/2$ for all points in $\exp(C_i)$. Then it follows that
$\exp (C_i)\cap E_{s_i}=\g(s_0-s_i)$. To see this, note that since
the second fundamental forms of a surface $\S\subset M$ and
$\exp^{-1}(\S)$ agree at the origin, we may reduce to the
Euclidean case. For two surfaces $\S_1,\S_2\subset \RR^3$, which
are tangent at $0$ to the horizontal plane, and such that
$\inf_{|e|=1}(A^{\nu_1}(e,e))>\sup_{|e|=1} (A^{\nu_2}(e,e))$,
where $\nu_j$ is the unit tangent to $\S_j$. Then locally near
$0$, $\S_1$ lies above $\S_2$. To see this, we note that in a
neighborhood of $0$, $\S_j$ is a graph over the horizontal plane.
The Hessian of $\S_j$ agrees with the second fundamental form at
the origin. Thus, the difference of the two graphs $\S_1-\S_2$ has
Hessian $>0$ in a neighborhood of the origin, and therefore has a
local max at the origin, which means that $\S_1$ lies above $\S_2$
(thanks to Ken Bromberg for help with this argument).

Now, we take the cap $\exp(C_i)$, and flow it along normal
geodesics distance $s_i$ to a surface $C_i'$. Since $\exp(C_i)\cap
E_{s_i}=\g(s_0-s_i)$, we see that $C_i'\cap E_0=\g(s_0)$. By the
equation $|Y'(s_i)| \leq  -\k s_{\k}(s_i) + |Y'(0)|$ (where $Y(t)
$ is now the Jacobi field associated to the exponential map from
$\exp(C_i)$), we see that if $s_i$ is small enough, then $C_i'$
will remain convex, since $Y'(0)$ cannot be zero. In fact, flowing
a bit further, we get a convex cap intersecting $E_0$ in a compact
set. We may then cut off a bit of $E_0$ by this cap to get a
smaller convex core, a contradiction to the minimality of $\CH(M)$
(the argument here is actually very similar to that of lemma
\ref{isoperimetric}). Thus, $\inf \{|W_{s_0s}e|\ \mid e\in T_p
C_{s_0}, |e|=1\}$ goes to zero as $s\to 0$. Since
$GK_{C_s}(r_{s_0s}(p)) Jac(r_{s_0s})(p)=Jac(-W_{s_0s})$, claims 1
and 2 follow from the bounds on $W_{s_0s}$.

 One can
thus make sense of the equation $K=K_{sec} + GK_s$ at points where
$C_s$ is $C^2$, where $K$ is the intrinsic sectional curvature of
$C_s$, and $K_{sec}$ is the sectional curvature of the plane
tangent to $C_s$. Since $K_{sec} \leq b<0$, we have
$$ 2\pi \chi(\d E_s)=
\lim_{s\to 0^+} \int_{C_s} K da = \lim_{s\to 0^+} \int_{C_s}
K_{sec}+ GK_{C_s} da $$ $$\leq \lim_{s\to 0^+} \int_{C_s} b da = b
\lim_{s\to 0^+} \area(C_s),$$ thus $ \area(C_0) \leq
\frac{2\pi}{b}\chi(\d \CH(M))$.

\end{proof}

\begin{theorem} ($4.5\frac34_+$ \cite{Gr99}) \label{systole}
Every closed, orientable surface $V$ of genus $g\geq 1$ with a
Riemannian metric  admits a closed curve of length $\leq l$ which
is not null-homologous and satisfies $\frac12 l^2\leq \area(V)$.
\end{theorem}

Now, we apply theorem \ref{systole} to conclude that each
component of $\d \CH(M)$ has an embedded curve of bounded length,
depending only on the topological type of $\d \CH(M)$.

%\bibliographystyle{abbrv}
%\bibliography{../refs}
\def\cprime{$'$} \def\cprime{$'$}

\end{document}